\documentclass[11pt]{article}

\usepackage[
  main=english, 
  english, german, russian, czech, slovak 
]{babel}        
\usepackage{paratype}

\usepackage{makeidx}      
\makeindex                  
\usepackage{paralist} 
\usepackage{amsmath}  
\usepackage{amsthm}
\usepackage{amsfonts}
\usepackage{amssymb}
\usepackage{graphicx}
\usepackage[all]{xy}
\usepackage{url}      
\usepackage{markdown} 
\usepackage{tabularx} 
\usepackage{tabu}
\usepackage{booktabs}
\usepackage{float}
\usepackage{multirow}
\usepackage{array}
\usepackage{cite}
\usepackage{tikz-cd}
\usepackage{csquotes}
\usetikzlibrary{arrows.meta}
\usetikzlibrary{cd}
\usepackage{titletoc}
\usepackage{CJK}
\usepackage{hyperref}
\usepackage{booktabs}
\usepackage{geometry}
\usepackage{listings} 
\newtheorem{theorem}{Theorem}[section] 
\newtheorem{proposition}[theorem]{Proposition}                              
\newtheorem{lemma}[theorem]{Lemma}         
\newtheorem{corollary}[theorem]{Corollary} 
\theoremstyle{definition}
\newtheorem{definition}{Definition}

\theoremstyle{remark}
\newtheorem*{remark}{Remark}

\begin{document}
\title{Explicit construction of Atiyah-Singer indices for maximally hypoelliptic operators on contact manifolds}
\author{Minjie Tian\\
 Department of Mathematics, Kyoto University}
\date{June 2021}
\maketitle
\begin{abstract}
The Atiyah-Singer index theorem gives a topological formula for the index of an elliptic differential operator. Enlightening from Alain Connes’ tangent groupoid proof of the index theorem and van Erp's research for the Heisenberg index theory on contact manifolds, we give an explicit construction of a series of maps, whose induced map in $K$-theory is the Heisenberg Atiyah-Singer index map on contact manifolds.

Our methods derive from Higson's construction for symbol class in $K$-theory.
\end{abstract}
\newpage
\tableofcontents
\titlecontents*{chapter}
[0pt]
{\addvspace{1em}}
{\bfseries\chaptername\ \thecontentslabel\quad}
{}
{\bfseries\hfill\contentspage}
\newpage
\section{Introduction}
Let $M$ be a compact oriented manifold.
Consider an elliptic operator $P$ between vector bundles $E$ and $F$ over $M$.
Adopting another compact space $X$, Atiyah and Singer\cite{atiyah1968index}
\cite{atiyah1968index2}  constructed 
two maps from $K(X\times TM)$ to $K(X)$ as the topological and analytical indices of $P$, respectively. They verified that these two indices coincide. When $X$ is a single point, the topological index is given by the integer $\dim\text{Ker} P-\dim \text{Coker} P$. The Atiyah–Singer index theorem is represented as
$$\int_{T^{\ast}M} Ch([\sigma(P)])\wedge Td(M)=\dim\text{Ker}(P)-\dim\text{Coker}(P).$$

This study aims to determine the images of such an analytical index map. In \cite{higson1993k}, Higson constructed  an asympototic morphism  $T^t$ from $C_0(T^{\ast}M)$ to $\mathcal{K}(L^2(M))$. For each $\alpha(x,y)\in C_c(T^{\ast}U)$, for some open subset $U$ of $M$ with $x,y$ ordinary coordinates, 
\begin{equation*}
\begin{aligned}
T^t_{\alpha}&:L^2(U)\to L^2(U),\\
T^t_{\alpha}f&(x)=\int_{T^{\ast}_x U} \alpha(x,t^{-1}y)e^{ixy}\hat{f}(y)dy
\end{aligned}
\end{equation*}
can be extended to an entire map over $C_0(T^{\ast}M)$.
Such $T^t$ induces an index map on symbols $Ind_P:K(X \times T^{\ast}M) \to K(X).$ It helps  to directly determine what the index map represents on the level of $C_0$ functions.

In \cite{van2010atiyah}\cite{van2004atiyah}, van Erp presented a particular construction of Heisenberg tangent spaces on Heisenberg manifolds and developed the index theory of hypoelliptic operators over them. 
 
 Consider hypoelliptic and maximally hypoelliptic operators over a contact manifold $M$, with $H\subset TM$ a codimension 1 subbundle.
 A vector field $X\in \Gamma(TM)$ is defined to have order 1 if $X\in \Gamma (H)$; otherwise, it is defined to have order 2.
  The product of the vector fields, $\Pi_i X^i$, has order $\sum d(i)$, where $d(i)$ is the order of $X^i$.
   Taking a local $H$-frame near $m$ on $M$, $ \{X^1,X^2, \cdots,  X^n\}$ ,
  each differential operator $P$ with Heisenberg order $d$ can be locally represented as 
  $$P=\sum_{|\alpha|\leq d} a_{\alpha}X^{\alpha},$$
  where $\alpha=\alpha_1+\cdots+2\alpha_n.$
 
 For the operators between bundles $E$ and $F$, $a_{\alpha}$ is taken as an $End(E,F)$-valued function over $M$, and $P$ can be represented by the same formula as above.
\begin{definition}
   A differential operator $P$ is called hypoelliptic if for any distribution $u$ on $M$, $Pu$ is smooth on an open $U\subset M$ implies that $u$ is smooth on $U$.
\end{definition}
\begin{definition}
  Let $(M,H)$ be a compact contact manifold.
  A differential operator $P$ on $M$ of $H$-order $d$ is called maximally hypoelliptic if for every differential operator $A$ on $M$ of $H$-order$\leq d$, there exists some constant $C$ satisfying \enquote{a priori} estimate
  $$\lVert Au\rVert_{L^2}\leq C(\lVert Pu\rVert_{L^2}+\lVert u\rVert_{L^2})$$ for all smooth $u\in C^{\infty}(M)$.
\end{definition}
Note that maximal hypoellipticity implies the hypoellipticity of the operators, and elliptic operators satisfy the 'a priori' estimate above; thus, all elliptic operators are hypoelliptic.
Let $M$ be a contact $(2n+1)$-dimensional manifold.
In fact the isomorphisms
   $$K_0(C^{\ast}(TM))\xrightarrow{\cong}K_0(C^{\ast}(\mathbb{T}_H M^{\text{adb}}))\xleftarrow{\cong}  K_0(C^{\ast}(T_H M)).$$
 are given there. 
The composition of the isomorphisms is denoted as
 $\Psi:K_0(C^{\ast}(T_HM))\cong K_0(C^{\ast}(TM))$.
 
 The main result of this study is the following.
 Take an open subset $U\subset M$, diffeomorphic to $G=H_n$, the Heisenber algebra. Let $\alpha$ and $\tilde{\alpha}$ be two idempotent matrix-valued functions on $T_H M$ and $TM$, supported on $T_H U$ and $TU$ respectively.  Thus, we can represent the functions on the restriction on $U$ by $\alpha(x,y)$ and $\tilde{\alpha}(x,y)$ for $x,y\in G$, where 
 $T_H G\cong G\times G$ and
$TG \cong G\times \mathfrak{g}\cong G\times G$ are obtained using the diffeomorphisms.
\begin{theorem}
\label{thm2}
 If
$$\lvert\alpha(x,\delta_t y)-\tilde{\alpha}(x,t y)\rvert= O(t^{-k})$$
for some integer $k>0$ as $t\to\infty$, 
then their classes in the $K$-theory correspond to each other by the isomorphism,
 $$\Psi:K_0(C^{\ast}(T_HM))\cong K_0(C^{\ast}(TM)).$$
\end{theorem} 
In fact, for generator functions $\{e^{-|x|^2}, x_1 e^{-|x|^2},\dots,x_n e^{-|x|^2} \} $ for $C_0(T_H G)$, we can typically choose appropriate $\tilde{\alpha}$ because they collapse to zero exponentially. Moreover, for arbitrary idempotent matrix-valued $\alpha\in C_0(T_H M)$, with a partition of unity, we may also obtain the corresponding $\tilde{\alpha}\in C_0(TM)$ as above. Such a correspondence is denoted by a map $\sim:C^{\ast}(T_H M)\to C^{\ast}(T M)$ extending to $C^{\ast}$-algebras.

Using the isomorphism, $\Psi$, we construct a similar map $T^t_H$ for hypoelliptic operators over contact manifolds.
Taking $\alpha(x,y)$ as a smooth, compactly supported function on $T_H G$, let
\begin{equation*}
    \begin{aligned}
     &{T^t_H}_{,\alpha}:L^2(G)\rightarrow L^2(G)\\
     &{T^t_H}_{,\alpha}(f(x))=t^{2n+1}\int_G \alpha(x,\delta_t(xy^{-1}))f(y)dy.
    \end{aligned}
\end{equation*}
Extending it to the entire $C^{\ast}(T_H M)$ yields the following:
\begin{theorem}
For a maximally hypoelliptic operator $P$, there exists an asymptotic morphism
 \[T^t_H:C^{\ast}(T_H M)\to\mathcal{K}(L^2(M))
 \]
 from the $C^{\ast}$-algebra of the osculating group $T_H M$ to $\mathcal{K}(L^2(M))$, such that it induces the Heisenberg topological index in the $K$-theory.
  \end{theorem}
 Actually, the construction of $T^t_H$ involves the analytical index, and the subsequent proof in section 6 can also be viewed as an explicit proof of the Atiyah–Singer index theorem.
 
 Using the explicit representation of $T^t$, we obtain several simple results:
 
 For a maximally hypoelliptic operator $P$, choose some $\sigma_H(P)_0$ to be represent $[\sigma_H(P)]\in K_0(C^{\ast}(T_H M))$. Such $\sigma_H(P)_0$ shall be considered as the difference of some idempotent matrix-valued functions.
 By Theorem \ref{thm2}, the induced map $T^{\ast}$ maps the symbol class $[\sigma_H(P)]$ to the Atiyah-Singer index. So the image of $\sigma_H(P)_0$ by $T^t$ shall be some formal difference of the compact operators, and if we take the rank we get the Atiyah-Singer index. Thus,
   \begin{corollary}
 
 \item[(1)]$rank ({T^t_H}_{,\sigma_H(P)_0})\in \mathbb{Z},$

\item[(2)] $\frac{\partial}{\partial t}rank({T^t_H}_{,\sigma_H(P)_0})=0,$

\item[(3)]Moreover, if $\sigma_H(P)_0$ can be taken to be compactly supported on $T_H U$ for some $U\cong G$, with $\sigma_H(P)_0(x,\delta_t(xy^{-1})=o(t^{-2n-1})$, 
 then 
 $Ind(P)=0$.
 \end{corollary}
Because $T^t$ is constructed from the contact structure and the Heisenberg index map is always non-trivial, we have the following:
 \begin{corollary}
 For a Heisenberg manifold $M$, if the collection of the homotopy classes of the asymptotic morphisms $[C^{\ast}{(T_H M)}, \mathcal{K}(L^2(M))]=0$, then $M$ cannot admit a contact structure.
 \end{corollary}
Although $T_H M$ and $TM$ are isomorphic, they have different topological and differential structures.  For contact manifolds, whose parabolic tangent bundles have more natural structures than the ordinary tangent bundle, the calculation may be performed directly over the parabolic ones. For example,
  
    Consider the following (non-commutative) diagram involving $T^t_H$ and  $T^t$:
    \[
\begin{tikzcd}
    &C^{\ast}(T_H M) \arrow{r}{\sim} \arrow{d}{T^t_H}
    &C^{\ast}(TM)\arrow{ld}{T^t}\\
    &\mathcal{K}(L^2(M))
\end{tikzcd}
 \]
  \begin{theorem}
The above diagram induces a commutative diagram on the level of the $K$-theory as follows:

   \[
\begin{tikzcd}
    &K_0(C^{\ast}(T_H M)) \arrow{r}{\Psi} \arrow{d}{Ind_H}
    &K_0(C^{\ast}(TM))\arrow{ld}{Ind}\\
    &\mathbb{Z}=K_0(\mathcal{K})
\end{tikzcd}
 \]
     \end{theorem}
\begin{corollary}
The index formula for a maximally hypoelliptic operator $P$ over a contact manifold is
given by
$$\int_{T^{\ast}M} Ch(\Psi[\sigma_H(P)])\wedge Td(M)=\dim\text{Ker}(P)-\dim\text{Coker}(P).$$
\end{corollary}

 As shown in \cite{willett2020higher}, the image of the operators $T^t$ can be considered to be an element in the localization algebra of the space, $L^2(T M,E\oplus F)$, when acting on $\alpha(\sigma(D))$ for $\alpha\in C_0(\mathbb{R})$. Higson proved that
 $$T^t(\alpha(\sigma(P)))-\alpha(t^{-1}P)\to 0$$
 for an elliptic $P$ of order $1$. When $P$ is formally self-adjoint and $M$ is complete, the latter term, $\alpha(t^{-1}P)$, represents the element $[P]\in K_0(M)=K_0(L^{\ast}(L^2(M)))$. Thus, on extending to a higher order, the image of $T^t$ represents the $K$-homology class $[P]$.
 
The same approach may be adopted for the hypoelliptic case, however the construction of $[P]$ in \cite{willett2020higher} cannot be extended directly because of the different definition of $\sigma_H$ for hypoelliptic operators. 

As a further research direction, the constructed asymptotic morphisms may also be performed in the $E$-theory. We also expect to construct a map whose images lie in the localization algebra to represent the elements in $K$-homology or an $E$-theory element, such that it induces the higher index map.

In the first half of this paper, we present a review of the theory of contact manifolds and hypoelliptic operators and explain several results. In section 5, we provide the definition and properties of the parabolic tangent groupoid. Moreover, in the last section, we describe the construction of the $T^t_H$ map and its induced relation with the index maps.

 \newpage
\section{Basic index theorem} 
\subsection{Symbol and index theorem}

   Let $P$ be a differential operator between the smooth
sections on vector bundles $E$ and $F$ over a smooth $n$-manifold $M$, i.e.,
$$P : L^2(M,E) \longrightarrow L^2(M,F).$$

 Choose a local coordinates of $M$ and perform a local trivialization of $E$ and $F$ near $m \in M$; $P$ can be locally written as
  $$ P=\sum_{|\alpha|= d}a_{\alpha}(x)\partial^{\alpha}_x . $$
  
  To understand the properties of $P$ near $m$, consider coordinate $x=t y$ is taken, and let $t$ converge to $0$, i.e.,
$$t^d P=\sum_{j=0}^d t^{j}P_{j},$$ 
where
$$P_{j}=\sum_{|\alpha|= d-j}a_{\alpha}(t y)\partial^{\alpha}_y,$$ 
and when $t\rightarrow 0$, $t^d P$ converges to
$$P_m=\sum_{|\alpha|= d}a_{\alpha}(0)\partial^{\alpha}_y .$$
It can be viewed as a map on the tangent space $T_m M$, with a constant coefficient for each $m$.

  From the operator $P_m$, we may construct the symbol at $m$ as
  $$\sigma (m,\xi)=\sum_{|\alpha|= d}a_{\alpha}(0) \xi^{\alpha},$$
  where $\xi=(\xi_1,\xi_2,\cdots, \xi_n)$ represents the regular local coordinates for $T_m M$, and $\xi=\sum_i {\xi^i x_i}$.
  Fixing $\xi$, $\sigma (m,\xi)$ will be a map from $E_m$ to $F_m$, with the entire map $\sigma (m,\xi)$ being a bundle map,
  $$\sigma: \pi^{\ast} E\longrightarrow \pi^{\ast} F,$$
  where $\pi:T^{\ast}M \rightarrow M$ is the projection.
     $P$ is defined to be elliptic if $\sigma (m,\xi)\in End(E_m,F_m)$ is invertible whenever $\xi\neq 0$. Such a bundle map yields the homotopy class 
     $$[ \sigma(P) ]\in K^0(T^{\ast}M).$$
    The Atiyah–Singer index theorem is expressed as
    $$\text{Index}\ P=\int_{T^{\ast} M} Ch([ \sigma(P) ])\wedge Td(M),$$
where $Ch$ is the Chern character and $Td(M)$ is the Todd class.
\subsection{$C^{\ast}$-algebra of groupoid}
\begin{definition}
 A groupoid is a category in which all morphisms are invertible.\\
 More specifically, denote the object set $G^{(0)}$ as a subset of the set of morphisms $G$, with several maps
 $$r:G\rightarrow G^{(0)},$$
 $$s:G\rightarrow G^{(0)},$$
 $$m:G^{(2)}\rightarrow ,G$$
 $$e:G^{(0)}\rightarrow G,$$
 $$i:G\rightarrow G,$$
which are the range, source, multiplication, unit , and inverse maps, satisfying the axiom of category.
 Here,
 $$G^{(2)}=\{(\gamma_1,\gamma_2)\in G\times G|s(\gamma_1)=r(\gamma_2)\}$$
 is the set of composable morphisms.
\end{definition}
In this study, we mainly consider a groupoid with a smooth structure and tangent groupoids.

We begin with the basic definitions of topological and smooth groupoids.
\begin{definition}
  A topological groupoid is a groupoid in which $G$ and $G^{(0)}$ are topological spaces, where $G^{(2)}$ is closed in $G$, and the above maps $r,s,m,e,i$ are continuous, $r$ and $s$ being open.
  A smooth groupoid is a groupoid in which $G$ and $G^{(0)}$ are smooth manifolds, where $G^{(2)}$ is closed in $G$, and the above maps $r,s,m,e,i$ are smooth, $r$ and $s$ being submersions.
\end{definition}
 Note that  $r$ and $s$ are open is equivalent to the condition that subsets $G_x=\{\gamma\in G^{(0)}|s(\gamma)=x\}$, $G^y=\{\gamma\in G^{(0)}|r(\gamma)=y\}$ are closed in $G$. Moreover, that $r$ and $s$ are submersions is equivalent to the condition that $G_x$ and $G^y$ are smooth submanifolds of $G$.\\
 \begin{definition}
   A Haar system on a topological groupoid $G$ is a family of positive measures $\{\lambda_x\}$ for each $x\in G^{(0)}$, defined on $G_x$, such that\\
   \begin{compactenum}
     \item Right invariance :$\forall$ $\gamma \in G$, $\forall A\subset G_{r(\gamma)}$,$$\lambda_{s(\gamma)}(A\cdotp \gamma)=\lambda_{r(\gamma)}(A).$$
     \item Continuity: $\int_{G_x}f\ d\lambda _x$ is continuous with respect to $x\in G^{(0)}$ for arbitrary $f\in C_c(G)$. 
     \end{compactenum}
 \end{definition}
 When we consider a smooth groupoid, $\lambda_{x}$ is taken to be a smooth 1-density, with $\int f\ d\lambda_x$ smooth for any $f\in C_c(G)$. 
 
 In the following, we only consider the smooth case.
For arbitrary $f,g\in C_c(G)$, the convolution and involution is defined as
$$(f\ast g)(\gamma)=\int_{\gamma^\prime \in G_x}f(\gamma{\gamma^{\prime}}^{-1})g(\gamma^{\prime})d\lambda_x ,$$
$$f^{\ast}(\gamma)=\overline{f(\gamma^{-1})},$$
where $x=s(\gamma)$. Such operations give $C_c(G)$ the $\ast$-algebra structure, and its completion for the norm 
$$\lVert f\rVert_{max}=\sup_{\pi}{\lVert\pi(f)\rVert}$$
yields the groupoid $C^{\ast}$-algebra, $C^{\ast}(G)$. Here, $\pi$ denotes the involutive Hilbert space representations of 
$C_c(G)$.

Below are some examples of smooth groupoids and their $C^{\ast}$-algebra.
 \begin{compactenum}
 \item Pair groupoid. $G=M\times M$, where $G^{(0)}=M$. For arbitrary $(x,y)\in M\times M$, $r(x,y)=x$ and $s(x,y)=y$, with composition $m((x,y),(y,z))=(x,z)$. In fact, we can consider the elements of the groupoid as lines connecting two arbitrary points over $M$, and the composition acting like a composition of paths. \\
 In this case, $\lambda_x$ should be a fixed smooth 1-density on $M$. The $C^{\ast}$-algebra have a canonical isomorphism:
 $$C^{\ast}(M\times M)\cong \mathcal{K}(L^2(M)).$$
 \item Groups. Given a group $\Gamma$, where $G = \Gamma$ and $G^{(0)} = {e}$, the law of composition is the group law.
 
 When we take a differential structure on $G$, since the object set is a single point, the measure is simply a smooth 1-density over $G$. The $C^{\ast}$-algebra in this case is just the group $C^{\ast}$-algebra.
 \item Tangent space $TM$ of a smooth manifold can be considered as a smooth groupoid, with the addition of vectors as group multiplication. It is a union of the groups with respect to the points on the manifold.
 
 A smooth Haar system on the tangent bundle $G = TM$ is a family of smooth densities on fibers $G_x = T_x M$ for $x\in M$. Considering the Haar measure on $T_x M$ as a linear map, $$\lambda_x : |\Lambda_n|T_x M \to\mathbb{R}^{+},$$
 the family $\{\lambda_x\}$ can be identified with a single $1$-density $\lambda$ on $M$. The $C^{\ast}$-algebra is 
 $$C^{\ast}(TM)\cong C^0(T^{\ast}M).$$
 \end{compactenum}

\subsection{Tangent groupoid}

For a smooth manifold $M$, Connes\cite{connes1994noncommutative} constructed the tangent groupoid, denoted as $\mathbb{T}M$, and proved the index theorem with that. 
Explicitly, a tangent groupoid is simply 
\begin{equation*}
\begin{aligned}
    \mathbb{T}M=G&=G_1\cup G_2\\
     &=TM \cup M\times M \times (0,1]\\
     &=\cup_{x\in M} T_x M \cup \cup_{t\in (0,1]} M\times M\times \{t\}.
\end{aligned}
\end{equation*}
It is the union of a series of tangent spaces (as an additive group) and pair groupoids, with a clear groupoid structure.

The smooth structure on $G$ is induced by the isomorphism near $t=0$.
\begin{equation*}
\begin{aligned}
    \psi: &TM\times [0,1]\rightarrow G,\\
          &(x,X,t)\mapsto (x,exp_X(-t X),t)\ \text{for}\ t>0,\\
          &(x,X,0)\mapsto (x,X)\ \ \ \ \ \ \text{for}\ t=0.\\
\end{aligned}
\end{equation*}
Such a smooth structure is equivalent to 
$$(x_t,y_t,t)\rightarrow (x,X)$$
if and only if
$$x_t,y_t\rightarrow x, \frac{x_t-y_t}{t}\rightarrow X.$$

Consider the $C^{\ast}$-algebra of $G,G_1$, and $G_2$ respectively. Thus,
$$C^{\ast}(G_1)\cong C_0(T^{\ast}M),$$
$$C^{\ast}(G_2)\cong C_0((0.1])\otimes \mathcal{K}(L^2(M)) .$$
These three $C^{\ast}$-algebras satisfy the following exact sequence:
$$0\rightarrow C^{\ast}(G_2)\rightarrow C^{\ast}(G)\xrightarrow{\sigma} C^{\ast}(G_1)\rightarrow 0.$$
Clearly $C^{\ast}(G_2)$ is contractible; therefore, $\sigma$ induces an isomorphism in the $K$-theory.
The target analytical index map is
$$Ind_a= \rho_{\ast}\circ (\sigma_{\ast})^{-1}:K^0(T^{\ast}M)\rightarrow \mathbb{Z}\cong K_0(\mathcal{K}),$$
where  $\rho$ is the restriction from $G_2=M\times M \times (0,1]$ to $M\times M \times \{1\}$. 

Connes proved in \cite{connes1994noncommutative}, using the tangent groupoid, that the analytical index map is the same as the topological index. A similar result can also be proved on the Heisenberg manifold for maximally hypoelliptic operators, which will be subsequently presented in section 6 of this paper.
\newpage
\section{Heisenberg structure and osculating group}
\subsection{Parabolic arrow}

Let $M$ be an $n$-manifold having a Heisenberg structure, i.e., with a distribution $H\subset TM$ with $\text{codim}\ H=1$.
van Erp\cite{van2010atiyah}\cite{van2004atiyah} provided a particular construction an the osculating group $G_m$, which is a unique simply connected Lie group corresponding to Lie algebra $\mathfrak{g}_m=H_m\oplus N_m$ for each $m\in M$, where $N=TM/H$ is the quotient line bundle. 
\begin{definition}
  Let $m\in U\subset  M$, and choose local coordinate $\Phi_U=(x_1.\cdots,x_n):U\rightarrow \mathbb{R}^n$, such that the first $n-1$ vectors $\{\frac{\partial}{\partial x_i}\}_{i=1,\dots, n-1}$ span $H$. Such a coordinate is called an $H$-coordinate.
\end{definition}
With such a coordinate,
\begin{definition}
  Take arbitrarily two smooth curves on $M$,
  $$c_1,c_2:[-1,1]\rightarrow M,$$
  tangent to $H$ at $t=0$.
  Define $c_1$ as equivalent to $c_2$ if in the $H$-coordinates,
  $$c_1^{\prime}(0)-c_2^{\prime}(0)=0,c_1^{\prime\prime}(0)-c_2^{\prime\prime}(0)\in H.$$
  This equivalence relation yields the equivalent class of $c$,
  denoted by $[c]_H$, which is called a parabolic arrow at $c(0)$.\\
  The set of all parabolic arrows at point $m\in M$ is denoted by $T_H M_m$. Moreover, the union,
  $$T_H M=\cup_{m\in M} T_H M_m,$$
  is the Heisenberg tangent bundle corresponding to the tangent bundle in the ordinary case.
\end{definition}
Such an equivalence class is independent of the choice of $H$-coordinates. The choice of the word "parabolic" can be understood as follows:

In the $H$-coordinate at $c(0)$, $c(t)$ may be written as
$$c(t)=c^{\prime}(0)t+\frac{1}{2}c^{\prime\prime}(0)t^2+O(t^3).$$
It can be equivalent to 
$$\tilde{c}(t)=ht+nt^2,$$
where $h=c^{\prime}(0)\in H$, $n=[\frac{1}{2}c^{\prime\prime}(0)]\in TM/H$.
Specifically, each $c(t)$ can be represented as a pair $(h,n)\in \mathbb{R}^n$, and it is called the Taylor coordinates for $[c]_H$ at $m$ with respect to $H$.

Compared with ordinary vector arrows on $M$, clearly parabolic arrows involve a second derivative relating to $H$.

 With such coordinates, $T_H M_m$ has the structure of a smooth manifold. 
 \begin{lemma}(\cite{van2004atiyah})
 \label{lem4}
 The tangent space at $[0]\in T_H M_m$ satisfies
 $$T_0(T_H M_m)\cong H_m\oplus N_m$$
 for each $m\in M$.
 \label{1emma 1}
 \end{lemma}
 Similar to multiplication by real numbers, dilation can be defined on $T_H M_m$.
\begin{definition}
  Dilation on $T_H M_m$ is defined as 
  $$\delta_s([c]_H)=[c_s]_H$$
  for arbitrary $s>0$, where $c_s(t)=c(st)$.
\end{definition}
van Erp also showed that parabolic arrows at the same point are composable using the composition of the local flows; therefore $G_m=T_H M_m$ have a Lie group structure for each $m$. By Lemma \ref{lem4}, its Lie algebra is $\mathfrak{g}_m=H_m\oplus N_m$.

We consider the groupoid structure of $C^{\ast}(T_H M)$ similar to $C^{\ast}(TM)$. Groupoid $T_H M$ is the union of a family of additive groups $T_H M_x$. For the $C^{\ast}$-algebra,  each regular representation $\pi_x, x\in M$, where
$$\pi_x(f)\Phi=f \ast \Phi,f\in C_0^{\infty}(G),\Phi\in C_0^{\infty}(G_x),$$
yields the $\ast$-
homomorphism, $$\pi_x : C^{\infty}_0 (T_H M)\to C^{\ast}(T_H M_x),$$
where $\lVert f\rVert_{C^{\ast}(T_H M)}=\sup_{x\in M}\lVert\pi_x(f)\rVert_{C^{\ast}(T_H M_x)}$. Moreover the function $x\to \lVert \pi_x(f)\rVert_{C^{\ast}(T_H M_x)}$ is continuous on $M$ for each $f\in C(T_H M)$,  i.e., we may consider the elements in $C^{\ast}(T_H M)$ as a continuous section of the bundle over $M$ with fibers $\{C^{\ast}(T_H M_x)\}_{x\in M}$.

\subsection{Exponential maps on $T_H M$}
\begin{definition}
  Let $(M,H)$ be a manifold with a Heisenberg structure.\\
  An exponential map
  $$exp:T_H M\rightarrow M$$
  is a smooth map such that for each $(m,v)\in T_H M$,  $c(t)=exp(m,\delta_t v)$ defines a parabolic arrow $[c]_H=v\in T_H M_m$.
\end{definition}
Locally, an exponential map can be constructed as follows:

Let $\Phi_m$ be a choice of $H$-coordinates for each $m$, and $\Psi_m$ be the Taylor coordinates for $T_H M_m$ with respect to $\Phi_m$.
$$\Phi_m:\mathbb{R}^n\rightarrow M,$$
$$\Psi_m:\mathbb{R}^n\rightarrow T_H M,$$
then $exp=\Phi_m\circ \Psi_m^{-1}$ is a Heisenberg exponential map.
Let $v\in T_H M_m$, which can be written in the Taylor coordinates, then 
$$\Psi_m^{-1}(v)=(h,n),$$
$$\exp(\delta_t(v))=\Phi_m(t h,t^2 n),$$
and the final result corresponds to the curve $c(t)=t h+t^2 h$ in the $H$-coordinates $\Phi_m$. Thus, 
$$[c]_H=(h,n)=v\in T_H M_m.$$

Another route to consider a Heisenberg exponential map is to consider the composition,
$$T_H M\xrightarrow{\log}H\oplus N\xrightarrow{j} TM \xrightarrow{\exp} M,$$
where $\log$ is the inverse of the Lie algebra exponential map
for all $m\in M$, and $j$ is induced by some choice of a section $N\hookrightarrow TM$.
Accordingly, each Heisenberg exponential map can be induced from the ordinary exponential map. 

\section{Contact manifolds and hypoellipticity}

\subsection{Differential operators on contact manifolds}
In this section, we consider manifolds with Heisenberg structures, and in particular, contact manifolds. The objective is to consider differential operators of particular order, called the Heisenberg order, to analyze the properties of hypoelliptic operators.

\begin{definition}
   Let $M$ be a smooth manifold of dimension $n$. A Heisenberg structure on $M$ is a subbundle $H\in TM$, where $\text{Codim}\ H=1$.
 
\end{definition}
  
 A vector field $X\in \Gamma(TM)$ is defined to have order 1 if $X\in \Gamma (H)$; otherwise it is defined to have order 2.
  A product of the vector fields, $\Pi_i X^i$, has order $\sum d(i)$, where $d(i)$ is the order of $X^i$.
  
  Taking a local $H$-frame near $m$ on $M$, $ \{X^1,X^2, \cdots,  X^n\}$ ,
  each differential operator $P$ with Heisenberg order $d$ can be locally represented as 
  $$P=\sum_{|\alpha|\leq d} a_{\alpha}X^{\alpha}.$$
  Here, $\alpha=\alpha_1+\cdots+2\alpha_n.$
  The span of all $\Pi_i X^i$ of order less or equal to $k$ is represented as $U^k$, with algebra $A^k=U^k/U^{k-1}$. Thus, $A=\oplus A^k$ can be viewed as the decomposition of $P$ with respect to the Heisenberg order.
  
  Take $P\in U^k$, let $P^k\in A^k$ denotes its principal part,  and more locally, consider the germ of $P^k$ near $m\in M$, i.e., $P^k_m\in\tilde{U}^k_m=A^k/\mathcal{I}_m A^k$, where
  $\mathcal{I}_m=\{f\in C^{\infty}(M)|f(m)=0\}$. The graded structure on $A$ can induce a graded structure on $\tilde{U}_m=\oplus \tilde{U}^k_m$.
  
To investigate the structure of graded algebra $\tilde{U}_m$, the following proposition is needed.

\begin{proposition}
    Take $\{X_i\}_{i=1.\cdots,n}$ as a local H-frame near $m\in M$.
    Let $Y_i=X_i(m)\in H_m$ for $i=1,\cdots,n-1$, $Y_n=X_n(m)^N\in N_m$; therefore $\{Y_i\}_{i=1,\cdots,n}$ forms a basis for $\mathfrak{g}_m=H_m\oplus N_m$. 

    Let $U(\mathfrak{g}_m)$ denote the universal enveloping algebra of $\mathfrak{g}_m$. Thus, the map,
    $$U(\mathfrak{g}_m)\rightarrow \tilde{U}_m,$$
    $$\sum a_{\alpha} Y^{\alpha}\mapsto \sum a_{\alpha}(X^{\alpha})^{|\alpha|}_m,$$
is an isomorphism of the graded algebra, and is independent of the choice of the $H$-frame.
\end{proposition}
    
    Locally, for arbitrary differential operator  $P=\sum_{|\alpha|\leq d} a_{\alpha} X^{\alpha}$, using the above proposition, $P_m=\sum_{|\alpha|=d} a_{\alpha}(m) Y^{\alpha}\in U(\mathfrak{g}_m)$ as above.
    Here, $\mathfrak{g}_m$ can be considered as the right invariant vector fields on $G_m=T_H M_m$; thus, $U(\mathfrak{g}_m)$ can be viewed as the right invariant differential operators on $G_m$.
    
    Therefore, for each $P$, a family of right invariant differential operators $\{P_m\}$ on $G_m$ is obtained, which are called the model operators for $P$ for each $m\in M$.

\subsection{Contact structure}
\begin{definition}
   Let $M$ be a manifold of dimension $2n+1$, where $H\subset TM$ is a codimension 1 subbundle. (Such an $M$ is called a Heisenberg manifold).
 
   $H$ defines a contact structure on $M$ if an arbitrary non-vanishing local 1-form $\theta$ on $M$ with $\theta(H)=0$ has the property that $\theta\wedge(d\theta)^n$ is a nowhere vanishing volume form.
\end{definition}

A canonical example of a contact manifold is the Heisenberg group, $H_n=\mathbb{R}^{2n+1}$, on which the group operator is given by
$$(x_0,\cdots,x_{2n})(y_0,\cdots,y_{2n})=(x_0+y_0+\frac{1}{2}(\sum_{j=1}^r(x_j y_{r+j}-y_j x_{r+j})),x_1+y_1\cdots,x_{2n}+y_{2n}).$$
On $H_n$, the right invariant vector fields are
\begin{equation*}
    \begin{aligned}
        X_0&=\frac{\partial}{\partial x_0},\\
        X_i&=\frac{\partial}{\partial x_i}+\frac{1}{2}x_{i+r}\frac{\partial}{\partial x_0},i=1,\cdots, n,\\
        X_{i+r}&=\frac{\partial}{\partial x_{i+r}}-\frac{1}{2}x_{i}\frac{\partial}{\partial x_0},i=1,\cdots, n,
    \end{aligned}
\end{equation*}
$H$ is generated by $\{X_i\}_{i=1,\cdots, n}$,
and the local 1-form $\theta$ is 
$$\theta=d x_0+\frac{1}{2}\sum_{i=1}^n(x_i d x_{i+n}-x_{i+n}d x_{i}),$$
which is called the canonical contact form on $H_n$.
The most important property of contact manifolds is the following, which is induced from Darboux's theorem.
\begin{proposition}
    Every contact manifold $(M, H)$ of dimension $2n+1$ is locally isomorphic (as a contact manifold) to an open subset of the Heisenberg group $H_n$ with its canonical contact structure, i.e. for contact form $\theta$ near $m\in M$, there exist local coordinates $\Phi=(x_0,\cdots,x_{2n})$, such that
    $$\theta=dx_0+\frac{1}{2}\sum_{i=1}^n(x_i d x_{i+n}-x_{i+n}d x_{i}).$$
\end{proposition}

\subsection{Rockland operators, hypoellipticity, and maximal hypoellipticity}
\begin{definition}
   A differential operator $P$ is called hypoelliptic if for any distribution $u$ on $M$, $Pu$ is smooth on an open $U\subset M$ implies that $u$ is smooth on $U$.
\end{definition}
This property is an easy consequence of the "a priori" estimate of the operator, which is used to define the following:

\begin{definition}
  Let $(M,H)$ be a compact contact manifold.
  A differential operator $P$ on $M$ of $H$-order $d$ is called maximally hypoelliptic if for every differential operator $A$ on $M$ of $H$-order$\leq d$, there exists some constant $C$ such that
  $$\lVert Au\rVert_{L^2}\leq C(\lVert Pu\rVert_{L^2}+\lVert u\rVert_{L^2})$$ for all smooth $u\in C^{\infty}(M)$.
\end{definition}
 When studying the local behavior of $P$, it may be approximated as $P_m$, which are homogeneous, right invariant differential operators on the graded nilpotent groups, $G_m=T_H M_m$.
 \begin{definition}
   A Rockland operator on a graded group $G$ is a differential operator $P$ satisfying the following:
\begin{compactenum}
      \item $P$ is right invariant and homogeneous,
      \item $d\pi(P)$ is injective on $S_{\pi}$ for each irreducible unitary representation on $G$, except for the trivial one.
    \end{compactenum}
    $S_{\pi}\subset H_{\pi}$ is the space of vectors $v\in H_{\pi}$ such that $g\mapsto \pi(g)v$ is a smooth map. Moreover, for representation $\pi$ on $G$, $d\pi$ is a representation on $\mathfrak{g}$ defined by
    $$d\pi(X)v=\frac{d}{d t}\Big|_{t=0}exp(tX)v$$
    for $X\in \mathfrak{g}$, $v\in S_{\pi}$.
 \end{definition}
 First, we present the properties of Rockland operators.
 \begin{theorem}(\cite{fischer2016quantization})
 Let $\mathcal{R}$ be a right-invariant, homogeneous differential operator on a graded group $G$. The hypoellipticity of $\mathcal{R}$ is equivalent to $\mathcal{R}$ satisfying the Rockland condition. In this case, any operator of the form,
 $$P=\mathcal{R}+\sum_{|\alpha|<d}c_{\alpha}X^{\alpha},$$
 where d is the Heisenberg order of $\mathcal{R}$, is also hypoelliptic.
 \end{theorem}
 \begin{corollary}
 Let $G$ be a graded Lie group, and
 $\mathcal{R}$ is the Rockland operator on $G$ of $H$-order $d$.
 There exist some $c>0$, such that 
 $$\forall\ \Phi\in \mathcal{S}(G), \lVert\Phi\rVert^2_{W^d_H}=\sum_{|\alpha|\leq d}\lVert X^{\alpha}\Phi\rVert_{L^2}\leq c(\lVert\mathcal{R}\Phi\rVert_{L^2}+\lVert\Phi\rVert_{L^2}).$$
 \end{corollary}

The relation between maximally hypoelliptic operators and Rockland operators over contact manifolds is the following:
\begin{theorem}(\cite{van2004atiyah})
  Let $(M,H)$ be a compact contact manifold.\\
  A differential operator $P$ on $M$ is maximally hypoelliptic if and only if all model operators $P_m,m\in M$ are Rockland operators.
\end{theorem}

\subsection{Symbol of ordinary elliptic operator in $K$-theory}
Let $M$ be a smooth, closed manifold, with an elliptic differential operator
$$P:C^{\infty}(M,E)\to C^{\infty}(M,F),$$
 with symbol $\sigma:\pi^{\ast}(E)\to\pi^{\ast}(F)$, which is a bundle map.\\
If we consider the exact sequence of the bundle map over $T^{\ast}M$, i.e.,
$$0\to \pi^{\ast}(E)\to\pi^{\ast}(F) \to 0,$$
then it yields the $K$-theory class, $[\sigma]\in K^0(T^{\ast}M)$,\\
Let $\bar{\sigma}=
\begin{pmatrix}
0 & -i \sigma\\
i\sigma^{\ast} & 0\\
\end{pmatrix}
$, $u=(\bar{\sigma}+i)(\bar{\sigma}-i)^{-1}$, $
\epsilon=
\begin{pmatrix}
-1 & 0\\
0 & 1
\end{pmatrix}$
$$[\sigma]=[\frac{1}{2}(\epsilon u+1)]-[\frac{1}{2}(\epsilon+1)]\in K^0(T^{\ast}M)$$
is the class above.

If we conduct the same construction for
$D=\begin{pmatrix}
0 & -i P\\
i P^{\ast} & 0\\
\end{pmatrix}
$, $U=(D+i)(D-i)^{-1}
$, $
\epsilon=\begin{pmatrix}
-1 & 0\\
0 & 1
\end{pmatrix}
$.
Thus, 
\begin{equation*}
    \begin{aligned}
    &[\frac{1}{2}(\epsilon U+1)]-[\frac{1}{2}(\epsilon+1)]\\
    =& [\ker D] -[\ker D^{\ast}]\in K_0(\mathcal{K}(H))
    \end{aligned}
\end{equation*}
is the Atiyah-Singer topological index of $D$.

\subsection{Symbol of maximally hypoelliptic operator}
  For a maximally hypoelliptic operator $P$, the Heisenberg symbol, $\sigma_H(P)$, can be defined as $\{P_m\}_{m\in M}$, which is a family of Rockland model operators on $G_m$.
   
Let 
\begin{equation*}
D_m=\left(
\begin{array}{cc}
  0  &  -i P_m\\
  i P_m^{\ast} & 0
  \end{array}
  \right)
  :=\Gamma(\pi^{\ast}E\oplus \pi^{\ast}F)\rightarrow \Gamma(\pi^{\ast}E\oplus \pi^{\ast}F)
\end{equation*}
where $\pi:T_H M\rightarrow M$ is the projection, and $\pi^{\ast}E\oplus \pi^{\ast}F$ is the bundle over $G_m=T_H M_m$.

Formally, such a $\{D_m\}_{m\in M}$ family can be viewed as some operator on sections of the entire bundle over $T_H M$, and it varies continuously when $m$ moves on $M$.

For the Heisenberg case, $P_m$ is not a constant coefficient operator as in the ordinary case. Thus, different from the ordinary symbol map, the Heisenberg symbol is not a bundle map but acts on sections of a bundle, 

Take $u_m=(D_m+i)(D_m-i)^{-1}$, and the Cayley transform of $D_m$, $\epsilon=
\begin{pmatrix}
       -1  & 0 \\
       0  &  1
\end{pmatrix}
\in End(\pi^{\ast}E\oplus \pi^{\ast}F)$.
The Heisenberg symbol class is defined as
$$[\sigma_H(P)]=[\frac{1}{2}(\epsilon u+1)]-[\frac{1}{2}(\epsilon+1)]\in K_0(C^{\ast}(T_H M)).$$
The reason the class lies in $K_0(C^{\ast}(T_H M)$ is a consequence of Lemma \ref{lem1}.

When $H=TM$ and $P$ is elliptic, then $[\sigma_H(P)]$ is equivalent to the topological class associated with the principal symbol,
$$[(\sigma(P)),\pi^{\ast}E,\pi^{\ast}F]\in K^0(T^{\ast} M).$$
\newpage
\section{Parabolic tangent groupoid}

\subsection{Basic smooth structure on $T_H M$}
 Similar to Connes' construction of the tangent groupoid  $\mathbb{T}M=TM\cup M\times M\times (0,1]$, the parabolic tangent groupoid can be constructed.
   Let $\mathbb{T}_H M=T_H M\cup M\times M\times (0,1]$. It can be divided into a family of groupoids
   $$T_H M=\cup_{m\in M} T_H M_m=\cup_{m\in M} G_m ,$$
   $$M\times M\times (0,1]=\cup_{t\in (0,1]} M\times M\times \{t\}=\cup_{t\in (0,1]} G_t.$$
   
   If an exponential map $\exp:T_H M\rightarrow M$ is given, the smooth structure of $\mathbb{T}_H M$ is defined to be induced by the isomorphism,
   \begin{equation*}
   \begin{aligned}
       \Phi&: T_H M\times [0,1] \rightarrow \mathbb{T}_H M,\\
       \Phi&(m,v,t)=(\exp_m(\delta_t v),m,t)\in G_t,t>0,\\
       \Phi&(m,v,0)=(m,v)\in T_H M_m .
         \end{aligned}
   \end{equation*}
 Equivalently, such smooth structure can be realized by combining $T_H M$ and $M\times M\times (0,1]$ by convergence as follows:
 
 Let $a(t),b(t)$ be smooth curves in $M$, where initial point $a(0)=b(0)=m$, $a^{\prime}(0),b^{\prime}(0)\in H_m$.
 $$\lim_{t\rightarrow 0}(a(0),b(0),t)=[a]_H\ast [b]_H^{-1}\in T_H M_m,$$
where the $\ast$ on the right hand side denotes the Lie group multiplication on $T_H M_m$.
\subsection{Parabolic tangent groupoid for contact manifold}
 To understand the smooth structure of  $\mathbb{T}_H M$ more directly, it is advantageous to consider a contact manifold. Locally it can be simply considered as the Heisenberg group $G=H_n$. with the standard contact structure.
  
  First, by right translation, $TG$ can be trivialized as
  $$TG\cong G\times T_0 G\cong G\times \mathfrak{g}.$$
 For $G=H_n$, the Lie algebra $\mathfrak{g}$ is
 $$\{X_i,Y_i,Z|i=1,\cdots, n, [X_i,Y_i]=Z\}$$
 where $H\subset TG$ is a subbundle generated by $\{X_i,Y_i\}$, and $N=TG/H$ is chosen to be the span of $Z\in \mathfrak{g}$.\\
 Thus $G\times \mathfrak{g}\cong N\oplus H$. Because the Lie algebra of $T_H G_x$ is $H_x\oplus N_x$, the exponential map can be obtained by applying the Lie algebra to the Lie group, which is an identification because $T_H G_x$ is niopotent.
 \begin{equation*}
     \begin{aligned}
      TG\cong & G\times \mathfrak{g}\rightarrow T_H G=G\times G,\\
        \cong & H\oplus N,\\
              &(x,\xi)\mapsto (x,exp_x(\xi)).
     \end{aligned}
 \end{equation*}
 The exponential map $T_H G\rightarrow G$ is expressed as
    $$exp_x(y)=yx,$$
    where $x\in G$, and $y\in G\cong T_H G_x$ is a parabolic arrow at $x$.
    
Consequently, the smooth structure of $\mathbb{T}_H G$ can be induced from
\begin{equation*}
    \begin{aligned}
    &\Psi: T_H G\times [0,1]\rightarrow \mathbb{T}_H G,\\
    &\Psi(x,y,t)=(\delta_t(y)x,x,t)\in G\times G\times (0.1],t>0,\\
    &\Psi(x,y,0)=(x,y)\in T_H G_x.
    \end{aligned}
\end{equation*}
The equivalent convergence relation is
$(\gamma_1(s),\gamma_2(s),s)\in G\times G\times (0,1]$, which converges to $(p,g)\in T_H G_p$ when $t\rightarrow 0$ if and only if 
$$\lim_{s\rightarrow 0}\delta_s^{-1}(\gamma_1\gamma_2^{-1})=g,$$
$$\lim_{s\rightarrow 0} \gamma_1(s)=\lim_{s\rightarrow 0} \gamma_2(s)=p.$$

In fact, we have the following:
\begin{proposition}
    The parabolic tangent groupoid of $G=H_n$ is isomorphic to the semi-direct product
    $$T_H G\cong (G\times[0,1])\rtimes_{\alpha}G,$$
    where $G$ acts on $G\times [0,1]$ by $\alpha(y)(x,t)=(\delta=t(y)x,t).$
    The isomorphism is obtained by
    \begin{equation*}
        \begin{aligned}
        ((x,t),y)&\mapsto (\delta_t(x),x,t)\in G\times G\times (0,1],t>0,\\
         ((x,0),y)&\mapsto y\in G\cong T_H G_x .
        \end{aligned}
    \end{equation*}

\end{proposition}
    
 \newpage   
\section{Construction and main theorem}
\subsection{Introduction}
Let $M$ be a smooth contact manifold with dimension $m=2n+1$.\\
As in \cite{connes1994noncommutative}, if we denote $e_0$ and $e_1$ as the restriction maps of groupoid $\mathbb{T}M$ to the $t=0$ and $t=1$ parts, respectively, the composition, ${e_1}_{\ast}\circ {e_0}^{-1}_{\ast}$ induces the index map $ K^0(T^{\ast}M)\cong K_0(C_0(T^{\ast}M))\cong K_0(C^{\ast}(TM))\rightarrow \mathbb{Z}$, which diagramatic form is
 \[
 \begin{tikzcd}
&K_0(C^{\ast}(\mathbb{T}M))\arrow{d}{{e_1}_{\ast}} \arrow{r}{{e_0}_{\ast}} &K_0(C^{\ast}(TM))\arrow{r}{\cong} &K^0(T^{\ast}M) \\
     &K_0(C^{\ast}(M\times M))\arrow{r}{\cong} &K_0(\mathcal{K}(L^2(M))).
 \end{tikzcd}
 \]
Similarly, for a parabolic tangent groupoid, such a diagram becomes:
 \[
 \begin{tikzcd}
&K_0(C^{\ast}(\mathbb{T}_H M))\arrow{d}{{e_1}_{\ast}} \arrow{r}{{e_0}_{\ast}} &K_0(C^{\ast}(T_H M))\\
     &K_0(C^{\ast}(M\times M))\arrow{r}{\cong} &K_0(\mathcal{K}(L^2(M))).
 \end{tikzcd}
 \]
The Heisenberg analytical index is defined by $Ind_H={e_1}_{\ast}\circ {e_0}^{-1}_{\ast}$.

An explicit homomorphism of $\ast$-algebras $C^{\ast}(T_H M)\rightarrow \mathcal{K}(L^2(M))$ can be considered to be defined such that the induced map in the $K$-theory is the Heisenberg index map. Such a construction can help to understand the structures of the symbol maps and index map.

Indeed, a family of maps $T^t:C^{\ast}(T_H M)\rightarrow \mathcal{K}(L^2(M))$ can be constructed, such that as $t\to \infty$, $T^t$ acts like a $\ast$-homomorphism.

\begin{definition}
  Let $A,B$ be $C^{\ast}$-algebras.
  An asymptotic morphism from $A$ to $B$ is a family of functions $T^t:A\rightarrow B$, such that
  \item[(1)] $T^t(a)$ is continuous with respect to $a$ and $t$.
  \item[(2)] $\limsup_{t\to \infty}\lVert{T^t(a)}\rVert<\infty$ for each $a\in A$.
  \item[(3)] $\lim_{t\to \infty}\lVert T^t(a)+\lambda T^t(b)-T^t(a+\lambda b)\rVert=0$,\\
  $\lim_{t\to \infty}\lVert T^t(a^{\ast})-T^t(a)^{\ast}\rVert=0$,\\
  $\lim_{t\to \infty}\lVert T^t(a) T^t(b)-T^t(ab)\rVert=0$,
  and the convergence is uniform on the compact subsets of $A$.
\end{definition}
Such $T^t:A\to B$ induces a map in the $K$-theory as 
$$T^{\ast}:K(A)\to K(B).$$
In fact, for any projection $p\in A$, $T^t(p)$ will be a quasi-projection for a sufficient large $t$. Thus, we can find some projection $q^t$ in the same equivalent class with $T^t(p)$ in $K(B)$.
We define $T^{\ast}[p]=[q^1]$, and it can be easily proved that this definition is independent of the choice of projection $q^t$.

Our objective is to construct 
$$T^t:C^{\ast}(T_H M)\rightarrow \mathcal{K}(L^2(M)),$$
with the induced map
$$T^{\ast}:K(C^{\ast}(T_H M))\rightarrow K(\mathcal{K}(L^2(M)))\cong \mathbb{Z},$$
coinciding with the Heisenberg index map $Ind_H={e_1}_{\ast}\circ {e_0}_{\ast}^{-1}$.
If we consider the correspondence, $\Psi$, between $K(C^{\ast}(T_H M))$ and $K(C^{\ast}(TM))$, which is induced from the correspondence of the groupoid $C^{\ast}$-algebras, then
$$\sim:C^{\ast}(T_H M)\rightarrow C^{\ast}(TM).$$
By explicit calculation the diagram
\[
\begin{tikzcd}
    &K(C^{\ast}(T_H M)) \arrow{r}{\Psi} \arrow{d}{Ind_H}
    &K(C^{\ast}(TM))\arrow{ld}{Ind}\\
    &\mathbb{Z}
\end{tikzcd}
 \]
can be shown to be commutative. Thus, we may consider the $Ind_H$ map simply as the composition $Ind\circ \Psi$.\\
\subsection{Comparing with the ordinary index map}
The relation between the Heisenberg and ordinary symbol maps is the following:\\
Consider the combining of the groupoids, as in \cite{van2004atiyah},
\[\begin{tikzcd}[column sep=huge]
H\oplus N  
& TM \arrow[swap,dashed]{l}{\delta_t^{-1}} \\
T_H M \arrow[swap,dashed]{u}{s^{-1}}
& M\times M  \arrow[swap,dashed]{u}{s^{-1}} \arrow[swap,dashed]{l}{\delta_t^{-1}}
\end{tikzcd}
\]
which is a family of groupoids over $M\times [0,1]\times[0,1]$, called the adiabatic groupoid (\cite{nistor2003index}) of $\mathbb{T}_H M$. For simplicity, it can be understood as the union of groupoid $\mathbb{T}_H M\times (0,1] $ with the Lie algebroid $Lie(\mathbb{T}_H M)$, with the convergence relation as $(s,exp(s X))\to (0,X)$.\\
The induced relation on $C^{\ast}$-algebra is
\[
\begin{tikzcd}[column sep=huge]
K_0(C^{\ast}(TM))\arrow[swap,dashed]{r}\arrow{d}{\cong} &K_0(C^{\ast}(T_H M))\arrow{d}{\cong}\\
K_0(C^{\ast}(\mathbb{T}M))\arrow{d}{e_1^{\ast}}\arrow[hook]{rd} &K_0(C^{\ast}(\mathbb{T}_H M))\arrow{ld}{e_1^{\ast}}\arrow[hook]{d}\\
K_0(\mathcal{K}(L^2(M)) &K_0(C^{\ast}(\mathbb{T}_H M^{\text{adb}}))\arrow{r}{\cong}& K_0(C^{\ast}(H\oplus N)).
\end{tikzcd}
\]

van Erp proved that $K_0(C^{\ast}(\mathbb{T}M))\to K_0(C^{\ast}(\mathbb{T}_H M^{\text{adb}}))$, $K_0(C^{\ast}(\mathbb{T}_H M))\to K_0(C^{\ast}(\mathbb{T}_H M^{\text{adb}}))$, $ K_0(C^{\ast}(\mathbb{T}_H M^{\text{adb}}))\to K_0(C^{\ast}(H\oplus N))$, which are induced by the inclusion(or restriction) maps, are isomorphisms.\\
Note that 
$Ind_H$ is the composition, $$K_0(C^{\ast}(T_HM))\xrightarrow{\cong}K_0(C^{\ast}(\mathbb{T}_H M))\xrightarrow{e_1^{\ast}}K_0(\mathcal{K}(L^2(M)),$$
where $Ind$ is 
$$K_0(C^{\ast}(TM))\xrightarrow{\cong}K_0(C^{\ast}(\mathbb{T}M))\xrightarrow{e_1^{\ast}}K_0(\mathcal{K}(L^2(M)).$$
$e_1^{\ast}$( in the Heisenberg case, and similarly for ordinary one) can be viewed as the composition
$$K_0(C^{\ast}(\mathbb{T}_H M))\xrightarrow{\cong} K_0(C^{\ast}(\mathbb{T}_H M^{\text{adb}}))\to K_0(\mathcal{K}(L^2(M))$$
where the second arrow is induced by the restriction from the large groupoid  $\mathbb{T}_H M^{\text{adb}}$ to the $s=t=1$ part. The commutativity of the left diagram is equivalent the the right one. If we take
$K_0(C^{\ast}(TM))\xrightarrow{\Psi} K_0(C^{\ast}(T_H M))$ as the composition of the isomorphisms,
$$K_0(C^{\ast}(TM))\xrightarrow{\cong}K_0(C^{\ast}(\mathbb{T}_H M^{\text{adb}}))\xleftarrow{\cong}  K_0(C^{\ast}(T_H M)),$$
then the diagram on the right is clearly commutative because it involves only inclusion maps.

The correspondence  $K_0(C^{\ast}(TM))\xrightarrow{\Phi} K_0(C^{\ast}(T_H M))$ can be just induced by a correspondence 
$C^{\ast}(TM)\sim C^{\ast}(T_H M))$.

In fact, 
if we consider $T_H M$ and $TM$ locally as
$T_H G\cong G\times G$ and
$T G \cong G\times \mathfrak{g}\cong G\times G$, respectively, then the combining of the adiabatic groupoid, as described at the beginning of this section, will be 
$$G\times G\times (0,1]\times (0,1]\ni(\delta_t(y)x,x,1,t)\xrightarrow{t\to 0} (x,y)\in T_H G\cong G\times G,$$
$$G\times G\times (0,1]\times (0,1]\ni(s(y) x,x,s,1)\xrightarrow{s\to 0} (x,y)\in TG\cong G\times G.$$
Note that because we consider $x,y\in G$, the addition structure on $\mathfrak{g}\cong G$ will also represented by multiplication.

Therefore, if we denote idempotent matrix-valued functions as $\alpha\in C_0(T_H M)$ and $\tilde{\alpha}\in C_0(TM)$, which are supported on $T_H U$ and $TU$, respectively. With the above isomorphisms, the restriction of functions on $U$ are represented by $\alpha(x,y)$ and $\tilde{\alpha}(x,y)$ for $x,y\in G$. Clearly they correspond if and only if they can extend to a continuous function over the entire adiabatic groupoid, ${\mathbb{T}_H M}^{\text{adb}}$. 

In the local coordinates, the function over $G\times G\times (0,1]\times (0,1]$ is denoted as $F(x,y,s,t)$; thus the continuity of the entire function yields that for arbitrary $\epsilon>0$, there exists some $\eta>0$, such that when $s,t\leq \eta$,
$$\lvert\alpha(x,y)-F(\delta_t(y)x,x,1,t)\rvert< \epsilon,$$
$$\lvert\tilde{\alpha}(x,y)-F(s(y)x,x,s,1)\rvert< \epsilon.$$
Thus,
$$\lvert\alpha(x,\delta_t^{-1}y)-F(y x,x,1,t)\rvert< \epsilon,$$
$$\lvert\tilde{\alpha}(x,s^{-1}y)-F(y x,x,s,1)\rvert< \epsilon.$$
Since $F$ is continuous, $F(y x,x,1,t)$ and $F(y x,x,s,1)$
are homotopic as functions for $(x,y)$. If $\epsilon$ is sufficiently small, then $\alpha(x,\delta_t^{-1}y)$ and $\tilde{\alpha}(x,s^{-1}y)$ adopt the same $K$-theory class for sufficiently small $s$ and $t$.

Conversely, if we restrict $\alpha$ and $\tilde{\alpha}$ by $\lvert\alpha(x,\delta_t^{-1}y)-\tilde{\alpha}(x,t^{-1}y)\rvert\leq C_k t^{k}$ as $t\to 0$ for arbitrary integer $k>0$, then 
$$\lvert F(y x,x,1,t)-F(y x,x,t,1)\rvert<2\epsilon+C_k t^{k}$$
can be taken to be sufficiently small. Therefore, it is natural to extend them to the entire $F$ over $G\times G\times [0,1]\times [0,1]$. Thus, these $\alpha$ and $\tilde{\alpha}$ corresponds.

Note that $C^{\ast}$ algebra $C_0(R^n)$ is generated by functions $\{e^{-|x|^2},x_1 e^{-|x|^2},\dots, x_n e^{-|x|^2}\}$. When $|x|\to\infty$, these functions converge to $0$ exponentially, and thus, uniformly restricted by $C_k |x|^{-k}$ for some constant $C_k$ for each $k>0$. Therefore, we can actually choose the corresponding $\tilde{\alpha}_0\in C_0(TG)$ for each $\alpha\in C_0(T_H G)$. Such $\tilde{\alpha}_0$ can be taken to be a quasi-idempotent matrix-valued term, so an idempotent $\tilde{\alpha}$ can be induced, and it is restricted by $\tilde{\alpha}_0$.

Replacing $t$ by $t^{-1}$ , the above discussion can be concluded as follows:
\begin{theorem}\label{thm1}
Let $M$ be a contact manifold with dimension $2n+1$.
Consider an open subset $U\subset M$, diffeomorphic to $G=H_n$. Let $\alpha$ and $\tilde{\alpha}$ be two matrix-valued functions on $T_H M$ and $TM$, supported on $T_H U$ and $TU$, respectively.  Thus, we can represent the functions on the restriction on $U$ by $\alpha(x,y)$ and $\tilde{\alpha}(x,y)$ for $x,y\in G$, where we identify
 $T_H G\cong G\times G$ and
$TG \cong G\times \mathfrak{g}\cong G\times G$ through the diffeomorphisms.
 If
$$\lvert\alpha(x,\delta_t y)-\tilde{\alpha}(x,t y)\rvert= O(t^{-k})$$
for some integer $k>0$ as $t\to\infty$, 
then their classes in the $K$-theory correspond by the isomorphism,
 $$\Psi:K_0(C^{\ast}(T_HM))\cong K_0(C^{\ast}(TM)).$$
\end{theorem} 
\begin{remark}
For arbitrary $\alpha\in C_{0}(T_H M)$, we may take a partition of unity $\{\rho_i\}_{i\in I}$ with respect to $\{U_i\}_{i\in I}$, with each $U_i\cong G$. The corresponding function, $\tilde{\alpha}$, may be taken as $\tilde{\alpha}=\sum_{i\in I} \widetilde{\rho_i \alpha}$. The above isomorphism presents a correspondence between the $K$-theory classes of $\alpha$ and $\tilde{\alpha}$.

Note that in this section, all coordinates are represented by the Heisenberg group $G$, with its multiplication structure. However, in the calculation on $TM$ described in the subsequent sections, the addition structure of $TM$ will be adopted and for example the corresponding function of $\alpha(x,\delta_t(y z))$ will be $\tilde{\alpha}(x,t(y+z))$.

We fit an appropriate $k$ in the subsequent proof of the calculation.
\end{remark}
\subsection{Construction of the asymptotic morphism}
The diagram in section 6.1 shows that the map can be represented as the composition of
$$e_0^{-1}:C^{\ast}(T_H M)\rightarrow C^{\ast}(\mathbb{T}_H M),$$
which induces an isomorphism in the $K$-theory, as a smooth expansion from $T_H M$ to $\mathbb{T}_H M$, with
$$e_1:C^{\ast}(\mathbb{T}_H M)\rightarrow C^{\ast}(M\times M),$$
which is simply the restriction of the groupoid to the ${t=1}$ part.\\
  To simplify the problem, we first assume the case of Heisenberg group $G=H_n\cong \mathbb{R}^{2n+1}$, i.e., to construct the map locally, and subsequently combine them using the partition of unity.
  
  On the Heisenberg group $G_n$, simply denoted as $G$,
  as in section 5, $TG$ and $G\times G\times (0.1] $ are combined by a particular convergence
 \begin{equation*}
 \begin{aligned}
  &G\times G\times (0,1]\ni(x_t,y_t,t)  \xrightarrow {t\to 0} (p,q)\in G\times G\cong T_H G \\
\Longleftrightarrow & \delta_t^{-1}(x_t y_t^{-1})\to q ; x_t,y_t\to p.
 \end{aligned}
 \end{equation*} 
 Note that the term $\delta_t^{-1}(x_t y_t^{-1})$ corresponds to $\frac{x_t-y_t}{t}$ in the ordinary case.
 We take the $(x,{\delta_t(y)}^{-1}x,t)$ sequence in $G\times G\times (0,1]$, which converges to $(x,y)$ 
in $T_H G$ with the above topology on $\mathbb{T}_H G$.

To extend the functions on $C^{\ast}(T_H G)$ to $C^{\ast}(\mathbb{T}_H G)$, first take $\alpha(x,y)\in C_0(T_H G)\cong C_0(G\times G)$. Moreover, define $\widetilde{\alpha}(x,y,t)\in C^{\ast}(G\times G\times (0,1])$ to satisfy
$$\widetilde{\alpha}(x,{\delta_t(y)}^{-1}x,t)=\alpha(x,y)$$
for all $t>0$, i.e., 
$$\widetilde{\alpha}(x,y,t)=\alpha(x,\delta_t^{-1} (xy^{-1})).$$
Because the construction is for asymptoticity near infinity, the convergence is taken as $t\rightarrow \infty$ instead of $t\to 0$. Let $\alpha_0(x,y,t)=\widetilde{\alpha}(x,y,t^{-1})=\alpha(x,\delta_t (xy^{-1}))$.\\
Restricting to $t=1$, 
$$K(x,y)=\alpha(x,xy^{-1})$$
is just the kernel of the target operator in $C^{\ast}(G\times G)\cong \mathcal{K}(L^2(G)$.
However, such a map
\begin{equation*}
    \begin{aligned}
     C_0(T_H G)&\to \mathcal{K}(L^2(G)),\\
     \alpha(x,y)&\mapsto Op(\alpha(x,xy^{-1})),
    \end{aligned}
\end{equation*}
does not present a $\ast$-homomorphsim; therefore, it is reasonable to consider operators with kernel $\alpha_0(x,y,t)$ with some other transformation.

Consequently, the constructed maps are 
\begin{equation*}
    \begin{aligned}
     T^t:C_0(T_H G)&\to \mathcal{K}(L^2(G))\\
     \alpha(x,y)&\mapsto Op(t^{2n+1}\alpha(x,\delta_t(x y^{-1}))).
    \end{aligned}
\end{equation*}

Assume that this family of maps is actually an asymptotic morphism. When restricting to $t=1$, it is simply the composition, ${e_1}\circ {e_0}^{-1}$, for the above choice of $e_0$ and $e_1$. Therefore, its induced map in the $K$-theory is the analytical index map $T^{\ast}=Ind_H={e_1}_{\ast}\circ {e_0}^{-1}_{\ast}$.

Thus, it suffices to show that such maps present an asymptotic morphism and can be combined together to yield a global map, and that the induced map $T^{\ast}$ in the $K$-theory is consistent with the Atiyah–Singer topological index map.

Fixing $\alpha(x,y)$ as a smooth, compactly supported function on $T_H G$,  it is presented as a function on $G\times G$, yielding
\begin{equation*}
    \begin{aligned}
     &T^t_{\alpha}:L^2(G)\rightarrow L^2(G)\\
     &T^t_{\alpha}(f(x))=t^{2n+1}\int_G \alpha(x,\delta_t(xy^{-1}))f(y)dy.
    \end{aligned}
\end{equation*}
Each $T^t_{\alpha}$ is a compact operator since its kernel has a compact support.

\subsection{Properties of maps and their combination}

In this section, we present some properties of the $T^t$ family and describe the method to combine them to become operators over $L^2(M)$ for the entire manifold $M$.

\begin{proposition}
    \item[(1)] $T^t_{\alpha}$ are uniformly bounded.
    \item[(2)] If $\beta(x,y)$ is another function in $C_0(T_H G)$, $T^t_{\alpha} T^t_{\beta}-T^t_{{\alpha}\ast {\beta}}\rightarrow 0$ in the operator norm as $t\to \infty$. The $\ast$ denotes the multiplication in the convolution algebra of groupoids.
    \item[(3)] If ${\alpha}^{\ast}$ is the complex conjugate of ${\alpha}$, then $T^t_{{\alpha}^{\ast}}-(T^t_{\alpha})^{\ast}\rightarrow 0$ as $t\to \infty$.
\end{proposition}
\begin{proof}
\item[(1)]
Note that $K^t_{\alpha}=t^{2n+1}\alpha(x,\delta_t( x y^{-1}))$, and $x y^{-1}$ can be restricted when $|x-y|$ has an upper bound since $\alpha$ has a compact support. Therefore, $K^t_{\alpha}\leq constant. \frac{t^{2n+1}}{(1+t|x-y|)^N}$ for a sufficiently large constant.
By simple calculus, 
$$ (T^t_{\alpha}f,f)\leq  constant. \int_G |f(x)|^2 \int_G \frac{t^{2n+1}}{(1+t|x-y|)^N} d y d x\leq constant.  \lVert f\rVert_{L^2}^2.$$
Thus, $T^t_\alpha$ is uniformly bounded.
\item[(2)]
Similarly to (1), such lemma can be proved.
\begin{lemma}(\cite{higson1993k})
\label{lem3}
Suppose that $A^t:L^2(G) \to L^2(G)$ are operators of the form, 
$$A^t f(x) = \int_G K^t(x, y)f(x) dy,$$
where $|K^t(x, y)|\leq constant. \frac{t^{2n+1}}{(1+t|x-y|)^{2n+2}}$.
For $L>0, A^t_L$ be an operator with kernel
\begin{equation*}
K^t_L(x,y)=
\begin{cases}
K^t(x,y) & \text{if $|x-y|<\frac{L}{t}$}\\
$0$ & \text {if $|x-y|\geq\frac{L}{t}$}.
\end{cases}
\end{equation*}
\end{lemma}
Thus, $\lVert A^t \to A^t_L\rVert \to 0 $ as $L\to \infty$, uniformly in $t$.
\begin{lemma}(\cite{higson1993k})
Let $A^t$ be as above, but suppose that
$$|K^t(x, y)| \leq constant.\  \frac{t^{2n}}{(1 + t|y-x|)^{2n+2}}.$$
Thus, $\lVert A^t\rVert\to 0$ as $t\to \infty$.
\end{lemma}
Thus, for some operators $A^t$ with kernels of compact support having properties in Lemma \ref{lem3}, to prove that $A^t\to 0$, we should prove the kernels of $A^t$ are bounded by a multiple of $t^{2n}$ on the set $\{x,y\big\lvert |x-y|< \frac{L}{t}\}$.

We can calculate the kernel directly, and the $O(t^k)$ following are restricted for $t\to \infty$.
\begin{equation*}
    \begin{aligned}
    &(T^t_{\alpha} T^t_{\beta} f)(x)\\
    =&t^{{2n+1}}\int_G \alpha(x,\delta_t(x y^{-1}))(T^t_{\beta}f)(y)d y\\
    =&t^{2({2n+1})}\int_G \alpha(x,\delta_t(x y^{-1}))\int_G \beta(y,\delta_t(y z^{-1}))f(z)d z d y\\
    =&t^{2 ({2n+1})}\int_G \alpha(x,\delta_t(x z^{-1}))\int_G \beta(z,\delta_t(z y^{-1}))f(y)d z d y.
    \end{aligned}
\end{equation*}
Thus, the kernel of $T^t_{\alpha} T^t_{\beta}$ is 
\begin{equation*}
    \begin{aligned}
K^t_1(x,y)=&t^{2({2n+1})}\int_G \alpha(x,\delta_t(x z^{-1}))\beta(z,\delta_t(z y^{-1}))dz\\
=&t^{2({2n+1})}\int_G \alpha(x,\delta_t(x y^{-1}z^{-1}))\beta(z y,\delta_t(z))dz\\
=&t^{2({2n+1})}\int_G (\tilde{\alpha}(x,t(x-y-z))+O(t^{-2n-2}))(\tilde{\beta}(z +y,t z)+O(t^{-2n-2}))dz\\
=&t^{2({2n+1})}\int_G \tilde{\alpha}(x,t(x-y-z))\tilde{\beta}(z +y,t z)dz+O(t^{2n}).
\end{aligned}
\end{equation*}

Because the multiplication is a convolution in groupoid $C^{\ast}$-algebra,
$$\alpha\ast \beta(x,y)= \int_G a(x,y z^{-1})\beta(x,z)dz;$$
therefore, 
\begin{equation*}
    \begin{aligned}
    &T^t_{\alpha\ast\beta}f(x)\\
    =&t^{2n+1} \int_G (\alpha\ast\beta)(x,\delta_t(xy^{-1}))f(y)dy\\
    =&t^{2n+1} \int_G \int_G \alpha(x,\delta_t(xy^{-1})z^{-1})\beta(x,z)f(y)dzdy.
    \end{aligned}
\end{equation*}
Thus, the kernel of $T^t_{\alpha\ast\beta}$ is 
\begin{equation*}
\begin{aligned}
K^t_2(x,y)&=t^{2n+1} \int_G \alpha(x,\delta_t(x y^{-1})z^{-1})\beta(x,z)dz\\
&=t^{2n+1} \int_G (\tilde{\alpha}(x,t(x-y)-z)+O(t^{-2n-2}))(\tilde{\beta}(x,z)+O(t^{-2n-2})dz\\
&=t^{2n+1} \int_G \tilde{\alpha}(x,t(x-y)-z)\tilde{\beta}(x,z)dz+O(t^{-1}).
\end{aligned}
\end{equation*}
The difference of the terms is
\begin{equation*}
\begin{aligned}
 &\int_G t^{2n+1}\tilde{\alpha}(x,t(x-y-z))\tilde{\beta}(z+y,t z)d z- \int_G \tilde{\alpha}(x,t(x-y)-z)\tilde{\beta}(x,z) d z\\
 = &\int_G \tilde{\alpha}(x,t(x-y)-z)(\tilde{\beta}(z/t +y,z)- \tilde{\beta}(x,z)) d z\\
 \leq & constant.\cdot t^{-1}
\end{aligned}
\end{equation*}
when $|x-y|<\frac{L}{t}$. Thus, $|K^t_1(x,y)-K^t_2(x,y)|\leq constant.\cdot t^{2n}$ as $t\to\infty$ over $|x-y|<\frac{L}{t}$. 

Thus, $T^t_{\alpha}T^t_{\beta}-T^t_{\alpha\ast\beta}(x,y)\to 0$.
Part \item[(3)] can be proved similarly.
\end{proof}
By extending $C_0(T_H G)$ to $C^{\ast}(T_H G)$, we obtain the following
\begin{theorem}
There is an asymptotic morphism
$$T^t:C^{\ast}(T_H G)\to \mathcal{K}(L^2(G)),$$
such that $\lVert T^t(\alpha)-T^t_{\alpha}\rVert\to 0$ for all $\alpha\in C_0(T_H G)$ as $t\to \infty$.
\end{theorem}
In the following we discuss how to combine the operators.

To generalize such an operator to an entire contact manifold, we first need the following:
\begin{lemma}
\label{lem2}
Let $f$ be a smooth, compactly supported function on G,
$M_f:L^2(G)\to L^2(G)$ is the left multiplication by f, then
$$M_f T^t_{\alpha}-T^t_{\alpha} M_f\to 0$$
as $t\to \infty$ in the operator norm.
\end{lemma}
It can be noticed from the definition that the integral takes values almost on the diagonal of $G\times G$, involving only the local part of the manifold, which explains why we can combine the operators.
\begin{proof}
The kernel of $M_f T^t_{\alpha}-T^t_{\alpha} M_f$ is
\begin{equation*}
    \begin{aligned}
    K^t(x,y)=&(f(x)-f(y))k^t(x,y).
    \end{aligned}
\end{equation*}
As $t\to \infty$, when 
we restrict it by 
$|x-y|\leq \frac{L}{t}$, 
$|f(x)-f(y)|< constant. \frac{L}{t}$; therefore, $|K^t|$ is restricted by $t^{2n+1}\cdot constant. \frac{L}{t}=constant.L \cdot{t^{2n}}$.
Thus, $M_f T^t- T^t M_f\to 0$ when $t\to\infty$.
\end{proof}
\begin{lemma}
Let $U$ and $V$ be two open subsets of $G$.
$\Phi:U\rightarrow V$ is a diffeomorphism, which also preserves the contact structure.\\
Let $\widetilde{\Phi}:G\times G\supset T_H U\to T_H V\subset G\times G$\\
$$\widetilde{\Phi}(x,y)=(\Phi(x),\Phi(y)).$$
$U_{\Phi}:L^2(U)\to L^2(V)$ is the induced isomorphism on $L^2$ spaces,
then $T^t_{\alpha\circ \widetilde{\Phi}}=U_{\Phi} T^t_{\alpha} U_{\Phi}^{-1}$ for all $\alpha\in C_0(T_H V).$
\end{lemma}
\begin{proof}
The kernel of $\widetilde{\Phi}=U_{\Phi} T^t_{\alpha} U_{\Phi}^{-1}$ is
$$K_1=t^{2n+1}\alpha(\Phi(x),\delta_t(\Phi(x)\Phi(y)^{-1})).$$
Moreover, the kernel of $T^t_{\alpha\circ \widetilde{\Phi}}$ is
$$K_2=t^{2n+1}\alpha(\Phi(y),\Phi(\delta_t(x y^{-1}))).$$
Since $\Phi$ preserves the contact structure,$\Phi(x)\Phi(y)^{-1}=\Phi(x y^{-1}),\Phi(\delta_t(x))=\delta_t(\Phi(x))$;
thus, $K_1=K_2$.
\end{proof}
Lemma \ref{lem2} suggests that the $T^t_{\alpha}$ maps only depend on the information near the diagonal, $x=y$, as $t\to \infty$. Thus, as $t$ becomes sufficiently large, $T^t_{\alpha}$ becomes commutative with the multiplication with the function, up to a small difference term.

If we relate two areas of $G$ by a contact diffeomorphism, preserving the contact structure, then the maps, $T^t_{\alpha}$ and $T^t_{\alpha\circ \widetilde{\Phi}}$ only differ by a unitary action of $U_{\Phi}$.

Therefore, using the partition of unity, we obtain the map,
$$T^t:C^{\ast}(T_H M)\rightarrow \mathcal{K}(L^2(M)),$$
which maps $\alpha$ to $\sum_{\gamma} U_{\Phi_{\gamma}}^{-1}T^t(\rho_{\gamma}\cdot \alpha)\circ \tilde{\Phi}_{\gamma})U_{\Phi_{\gamma}}$ for some  partition of unity $\{\rho_{\gamma}\}$ and local contact trivialization $\Phi_{\gamma}$. Moreover, as $t\to \infty$, the resulting map is independent of the choice of partition up to a small difference. As a conclusion,
\begin{theorem}
Let $M$ be a smooth contact manifold without boundary, with a Haar measure. There is an asymptotic morphism
$$T^t:C^{\ast}(T_H M)\to \mathcal{K}(L^2(M)),$$
such that if $\Phi:U\to G$ is a diffeomorphism-preserving contact structure from an open set in $M$ to $G$, then
$$\lVert T^t(\alpha\circ \widetilde{\Phi})-U_{\Phi}T^t_{\alpha}U_{\Phi}^{-1}\rVert \to 0$$
as $t\to \infty$, for all $\alpha\in C_0(T_H G)$
\end{theorem}

\subsection{$T^{\ast}$ and index map}

First, some concepts are explained.
Let $P$ be a maximally hypoelliptic operator, where $P\in U^k\setminus U^{k-1}$ involves only the $k$-order part for simplicity.
$$P:C^{\infty}(M,E)\to C^{\infty}(M,F).$$
In this case, instead of the ordinary symbol $\sigma\in K^0(T^{\ast}M)$, we consider the Heisenberg symbol $\sigma_H\in K_0(C^{\ast}(T_H M))$.

Recall that $P$ induces a smooth family of Rockland operators $P_m$ on $G_m$.
Let \begin{equation*}
    D_m=\left(
    \begin{array}{cc}
       0  & -i P_m \\
       i P_m     &  0
    \end{array}
\right) 
,
 D=\left(
    \begin{array}{cc}
       0  & -i P \\
       i P     &  0
    \end{array}
\right) ,
\end{equation*}
then $D_m$ acts on sections over the bundle, $\pi^{\ast}E\oplus \pi^{\ast}F$, over $G_m$, where $\pi:T_H M\to M$ is the projection.

Consider a continuous family, $\sigma_H(D)=\{D_m\}$, such a family can be viewed formally as a map between the sections of $\pi^{\ast}E\oplus \pi^{\ast}F$ over the entire $T_H M$, which is continuous with respect to $m$. This is just the Heisenberg symbol of $D$.

With $u_m=(D_m+i)(D_m-i)^{-1}$, the symbol class of $P$, which is denoted as $[\sigma_H(P)]$, is defined as 
$$[\sigma_H(P)]=[\frac{1}{2}(\epsilon u+1)]-[\frac{1}{2}(\epsilon+1)]\in K_0(C^{\ast}(T_H M)),$$
where $u$ is family $\{u_m\}$ and $\epsilon=
\begin{pmatrix}
-1 & 0 \\
0 &  1
\end{pmatrix}
$
is an involution of the bundle $\pi^{\ast}E\oplus \pi^{\ast}F$.

Suppose $\alpha \in C_0(\mathbb{R})$, 
$$\alpha(t^{-k}D):L^2(M,E\oplus F)\to L^2(M,E\oplus F)$$ is defined by the spectral theory as follows:
$$\alpha(t^{-k}D)u_{\lambda_n}=\alpha(t^{-k}\lambda_n)u_{\lambda_n}.$$
Since $|\lambda_n|\to \infty$ as $n\to \infty$, such operators are compact operators.

When dealing with $\alpha(\sigma_H(D))$ for some $\alpha\in C_0(\mathbb{R})$, different from $\alpha(\sigma(D))$, it involves the differential part; therefore, it is difficult to represent it as a matrix-valued function over $T_H M$ directly.

However, we have the following lemma, which was proved by van Erp.
\begin{lemma} (\cite{van2010atiyah})
\label{lem1}
Let G be a graded niopotent group. If we identify $C^{\ast}(G)$ under the left regular representation on $L^2(G)$, i.e., under $Uh=u\ast h$ ,  $u\in C_0(G)$ and $U\in B(L^2(G))$ can be viewed as the same elements in $C^{\ast}(G)$.

For any Rockland self-adjoint operator $L$ on $G$, for  $\forall\ \alpha\in C_0(\mathbb{R})$,
$$\alpha(L)\in C^{\ast}(G).$$
\end{lemma}

In this case, $L$ is the matrix of operators $D_m$; therefore $$\alpha(D_m)\in C^{\ast}(G_m)\otimes_M C(End(E\oplus F)).$$
Thus,
$$\alpha(\sigma_H(D))\in C^{\ast}(T_H M)\otimes_M C(End(E\oplus F)).$$
We may consider $\alpha(\sigma_H(D))$ as the convergence of some matrix-valued $C_0$ series of functions on $T_H M$. Therefore, $T^t$ can be applied to it.\\

The key property of $T^t$ is as follows:
\begin{proposition}\label{prop1}
    Let $P$ be a maximally hypoelliptic operator with only the k-order part. $D=
    \begin{pmatrix}
    0 & -i P\\
    iP^{\ast} & 0
    \end{pmatrix}
    $, with Heisenberg symbol $\sigma_H(D)$, then
    $$T^t(\alpha(\sigma_H(D)))-\alpha(t^{-k}D)\to 0$$
    in the operator norm as $t\to \infty$ for each $\alpha\in C_0(\mathbb{R})$.
\end{proposition}

\begin{proof}
First, we relate $\alpha(\sigma_H(D))\in C^{\ast}(T_H M)\otimes_M C(End(E\oplus F))$ to some term in $C_0(T^\ast M)$; therefore we can relate $\sigma_H(D)$ and the operator $D$.
The map is constructed as the composition,
$$\wedge:C^{\ast}(T_H M)\rightarrow C^{\ast}(TM)\xrightarrow{Fourier} C_0(T^{\ast}M),$$
in which the first map is the correspondence in section 6.2, and the second map is the Fourier transform over each tangent space.

To understand the map more in detail, we considered some calculations locally.

In the local Taylor coordinates, multiplication on $T_H G_m$ is obtained by
$$(h,n)\ast(h'.n')=(h+h',n+n'+b(h,h'))$$
for some 2-form $b$. In this case of Heisenberg group $H^{2r+1}$, with $h=(h_1,\cdots, h_{2r})$, $b=
\begin{pmatrix}
0 & I_r \\
-I_r & 0
\end{pmatrix}
$.

For differential operator $P=\sum_{|\alpha|\leq d}a_{\alpha}{\partial}^{\alpha}$, it is recalled that the model operator at $m$ is
$P_m=\sum_{|\alpha|\leq d}a_{\alpha}Y^{\alpha}$, where $Y^i={\partial}^i(m)\in U(g_m)$ denotes the right invariant vector fields on $G_m$. In fact, by direct calculus, $Y^i=exp_{\ast x}\partial_i=\partial_i+b(\partial_i,x)$.

Thus, for Heisenberg symbol $\sigma_H(P)=\{P_m\}_{m\in M}$, the corresponding right invariant operator on $TM$ is $\tilde{P}_m=\sum_{|\alpha|\leq d}a_{\alpha}(m)\partial^{\alpha}$. When taking the Fourier transform on the acting space $C^{\ast}(TM)\to C_0(T^{\ast}M)$, the operator transforms to $\sigma=\sum_\alpha a_{\alpha}(m)\xi^{\alpha}$ on $T^{\ast}_m M$, with natural coordinate $\xi=\sum_i x_i \xi^i$. This operator is the ordinary symbol of $P$; therefore, such a transformation is adopted. 

However, this transform $\wedge$ is performed on the acting space of the differential operator; however, as in Lemma \ref{lem1}, we view the operators as particular functions, or their convergence.\\
For some $p$, matrix-valued $C_0$ function  on $T_H M$, with operator $Pf=p\ast f$, such correspondence follows:
\[
\begin{tikzcd}
    &C^{\ast}(T_H M)\arrow{r} &C^{\ast}(TM)\arrow{r} &C_0(T^{\ast}M)\\
    &p\arrow[mapsto]{r} &\tilde{p}\arrow[mapsto]{r} &\hat{p}
\end{tikzcd}
\]

Viewed as a function, the resulting function $\hat{p}$ in $C_0(T^{\ast}M)$ is simoly the Fourier transform of the transformation of $p$ to $C^{\ast}(TM)$.

As an operator, we transform on the acting function spaces, and obtain the operators as
$$P\mapsto \tilde{P}\mapsto \hat{P},$$
where $\hat{P}(\hat{f}):=\widehat{P(f)}$. However, note that
\begin{equation*}
\begin{aligned}
\hat{P}(\hat{f}):&=\widehat{P(f)}\\
&=\widehat{\tilde{p}\ast \tilde{f}}\\
&=\hat{p}\ast\hat{f}.
\end{aligned}    
\end{equation*}
Thus, the corresponding function of $\hat{P}$ is just $\hat{p}$.

For $\alpha (\sigma_H(D))\in C_0(T_H M)$, the corresponding operator on $C_0(T^{\ast} M)$ is denoted as $\widehat{\alpha(\sigma_H)}$.
Using a partition of unity, we may simplify the problem to $\alpha (\sigma_H(D))\in C_0(T_H G)$, and denote the function as $\alpha (\sigma_H(D))(x,y)$ for $x,y\in G$.

As convergences of matrix-valued functions, the Fourier transform is
$$\widetilde{\alpha(\sigma_H D)}(x,y)=\int_G \widehat{\alpha(\sigma_H)}(x,\xi)e^{-i y\xi }d\xi.$$

Based on the local calculation discussed above, $\widehat{\alpha(\sigma_H)}$ is the ordinary symbol of $\alpha(D)$; therefore, for $u$ supported on a small open set, we have the relation
$$(\alpha (D))u)(x)=\int_G e^{i x\xi}\widehat{\alpha({\sigma_H})}(x,\xi)\hat{u}(\xi)d \xi.$$
For $u$ supported on a small open set $U\subset G$, by calculation,
\begin{equation*}
    \begin{aligned}
    (T^{t}_{\alpha(\sigma_H(D))}u)(x)&=t^{m}\int_G \alpha(\sigma_H)(x,\delta_t(x y^{-1}))u(y)d y\\
    &=t^{m}\int_G (\widetilde{\alpha(\sigma_H)}(x,t(x-y))+O(t^{-m-1}))u(y)d y\\
    &=t^{m}\int_G\int_G \widehat{\alpha(\sigma_H)}(x,\xi)e^{i t (x-y) \xi}u(y)d\xi d y+O(t^{-1}).\\
    (Du)(x)&=\int_G e^{i x\xi}\widehat{\sigma_H}(x,\xi)\hat{u}(\xi)d \xi,\\
    (\alpha(t^{-k}D)u)(x)&=\int_G e^{i x\xi}\widehat{\alpha(t^{-k}{\sigma_H})}(x,\xi)\hat{u}(\xi)d \xi\\
    &=\int_G\int_G \widehat{\alpha(t^{-k}{\sigma_H})}(x,\xi) e^{i x\xi} e^{-i y\xi}u(y) d y d\xi\\
    &=\int_G\int_G \widehat{\alpha({\sigma_H}(x,\delta_t y))}(x,\xi) e^{i x\xi} e^{-i y\xi}u(y) d y d\xi.\\
    &=\int_G\int_G \widehat{(\widetilde{\alpha({\sigma_H}(x,t y))}+O(t^{-m-1}))}(x,\xi) e^{i x\xi} e^{-i y\xi}u(y) d y d\xi\\
    &=\int_G\int_G \widehat{\alpha({\sigma_H})}(x,t^{-1}\xi) e^{i x\xi} e^{-i y\xi}u(y) +O(t^{-m-1})d y d\xi\\
    &=t^{m}\int_G\int_G \widehat{\alpha({\sigma_H})}(x,\xi) e^{i t x\xi} e^{-i t y\xi}u(y) d y d\xi+O(t^{-m-1})\\
    &=(T^t_{\alpha(\sigma_H(D))}u)(x)+O(t^{-1})
    \end{aligned}
\end{equation*}
as $t\to\infty$.
Thus,
$$T^t(\alpha(\sigma_H)-\alpha(t^{-k}D)\to 0$$ as $t\to \infty$. As $t$ becomes larger, the integral area is taken to be sufficiently small; therefore the result holds for all functions $u$ over $M$.
\end{proof}
Using this property, we obtain the following:
\begin{theorem}
$$T^{\ast}([\sigma_H(P)])=Ind(P)=\dim\text{Ker} P-\dim\text{Coker} P.$$
\end{theorem}
\begin{proof}
For $0<t<\infty$, let
$$U_t=\frac{t^{-k}D+i}{t^{-k}D-i}$$
be the action of function $\alpha(x)=\frac{x+i}{x-i}$ on $t^{-k}D$. 

Thus, by proposition \ref{prop1}, $U^t-T^t(u)\to 0$ as $t\to \infty$; therefore,
\begin{equation*}
    \begin{aligned}
    T^{\ast}([\sigma_H(P)])&=T^{\ast}[\frac{1}{2}(\epsilon u+1)-\frac{1}{2}(\epsilon +1))]\\
    &=[\frac{1}{2}(\epsilon U_1+1)-\frac{1}{2}(\epsilon +1))]\in K(\mathcal{K}(L^2(M))).
    \end{aligned}
\end{equation*}

However, as $t\to 0$, $U_t=\frac{t^{-k}D+i}{t^{-k}D-i}$ converges to the operator maps, $f\in\ker D$ to $-f$, and $f$ in the orthogonal part to itself. If we denote $[F]$ as the projection to the subspace $F$, then
\begin{equation*}
    \begin{aligned}
    U_0&=-[\ker D]+(1-[\ker D])\\
    &=1-2[\ker D].
    \end{aligned}
\end{equation*}
By the definition of $D$, $[\ker D]=
\begin{pmatrix}
[\ker P] & \ \\
\  & [\ker P^{\ast}]
\end{pmatrix}
$; therefore,
\begin{equation*}
    \begin{aligned}
    \frac{1}{2}(\epsilon U_0+1)&=\frac{1}{2}(\epsilon-2\epsilon[\ker D]+1)\\
    &=
    \begin{pmatrix}
    [\ker P] & \ \\
\  & 1-[\ker P^{\ast}]
 \end{pmatrix},\\
  \frac{1}{2}(\epsilon +1)&= [1-[\ker P^{\ast}]]+[\ker P^{\ast}].
    \end{aligned}
\end{equation*}
Thus, by the homotopy invariance of the $K$-theory,
\begin{equation*}
    \begin{aligned}
    T^{\ast}([\sigma_H(P)])&=[\ker P]-[\ker P^{\ast}]\in K(\mathcal{K}(L^2(M)))\\
    &=Ind(P).
    \end{aligned}
\end{equation*}
\end{proof}

Finally, we may examine the commutative diagram directly on the $C^\ast$-algebra level.

First, $T^t_H$ (to avoid confusion we use $H$ for the Heisenberg case) maps $\alpha(x,y)\in C_0(T_H G)$ to the compact operator with kernel $t^m\alpha(x,\delta_t(xy^{-1}))$, and as shown by Higson\cite{higson1993k}, $T^t$ maps $\hat{\alpha}(x,y)\in C_0(T^{\ast} G)$ to the compact operator with kernel 
$$(\frac{t}{2\pi})^m\int_G \hat{\alpha}(x,\xi)e^{i t (x-y)\xi} d\xi
,$$ inducing the ordinary topological index. It is noted that in fact it is simply $t^m (\mathcal{F}^{-1}\hat{\alpha})(x,t(x-y))= t^m\tilde{\alpha}(x,t(x-y))$.

The above is extended to functions $\alpha\in C_0(T_H M)$ and $\tilde{\alpha}\in C_0(T M)$ using the partition of unity. With the correspondence in Theorem \ref{thm1} and the Fourier transform over each fiber,
$$\mathcal{F}:C_0(T_x M)\to C_0(T^{\ast}_x M),$$
the following diagram commutes.
\[
\begin{tikzcd}
    &K(C^{\ast}(T_H M)) \arrow{r}{\Psi} \arrow{d}{Ind_H}
    &K(C^{\ast}(TM))\arrow{ld}{Ind}\\
    &\mathbb{Z}
\end{tikzcd}
 \]
 
Thus, the Atiyah-Singer index theorem for maximally hypoelliptic operators on contact manifolds is represented as
\begin{equation*}
    \begin{aligned}
    &Ind\ P=\dim\text{Ker}(P)-\dim\text{Coker}(P)\\
    =&Ind_H([\sigma_H(P)])\\
    =&Ind(\Psi[\sigma_H(P)])\\
    =&\int_{T^{\ast}M} Ch(\Psi[\sigma_H(P)])\wedge Td(M).
    \end{aligned}
\end{equation*}

\newpage
\bibliographystyle{plain}
\bibliography{example}

\end{document}